\documentclass[11pt,english]{article}

\usepackage{graphicx,psfrag,epsfig,amsfonts,amssymb,amsmath,babel}

\newtheorem{Proposition}{Proposition}
\newtheorem{Lemma}[Proposition]{Lemma}
\newtheorem{Definition}[Proposition]{Definition}

\newtheorem{Remark} [Proposition]{Remark}

\newtheorem{Theorem}[Proposition]{Theorem}

\def\dist{\mbox{dist}}
\def\diam{\mbox{diam}}

\begin{document}
\title{CHAOS AND STABILITY IN A MODEL OF  INHIBITORY NEURONAL NETWORK.}
\author{Eleonora Catsigeras   
\thanks{Instituto de Matem\'{a}tica y Estad\'{\i}stica Prof. Rafael Laguardia
 Facultad de Ingenier\'\i a, Universidad de la Rep\'{u}blica.  Uruguay. }}
 
\date{November 14th., 2008}

\maketitle
\begin{abstract}

  We analyze the dynamics of a deterministic model of inhibitory neuronal networks proving that the discontinuities of the Poincar\'{e} map produce a never empty chaotic set, while its continuity pieces produce stable orbits.
 We classify the systems in three types: the almost everywhere (a.e.) chaotic,  the a.e. stable,  and the combined systems. The a.e. stable are periodic and chaos  appears as bifurcations. We prove that a.e. stable systems exhibit  limit cycles, attracting a.e. the orbits.

\bigskip


\end{abstract}

\noindent{\bf INTRODUCTION}

\noindent We obtain  rigorous mathematical results on the dynamics of a deterministic abstract discontinuous dynamical system, in a finite but large dimensional phase space. It comes from a non linear model of  inhibitory neuronal networks, without delays, composed by equally or different $2 \leq n < +\infty$ pacemaker neurons, evolving according to an autonomous differential equation in the inter-spike interval times, and interacting among them by synaptical instantaneous currents in the spiking instants.

\noindent A vectorial autonomous differential equation  governs the increasing potentials of the $n$ neurons as a function on time $t$, only during the inter-spike regime. On the other hand, the spiking regime holds when at least one neuron, say $i$, reaches a threshold level and gives a spike.  Due to the synaptic connections among the neurons, this spike produces sudden changes in the potentials  of the other neurons  $j \neq i$ and resets the potential of the neuron $i$.

\noindent The synapsis is assumed to be inhibitory, i.e. phase redeeming, meaning that the potentials of the receiving neurons suffer  negative changes in the spiking times. That is why the inhibitory synapsis is modeled by a matrix $\{H_{ij}\}_{i \neq j}, \, 1 \leq i,j \leq n$ of negative numbers $H_{ij} < 0$ that represent the instantaneous discontinuity jumps in the potentials of the neurons $j \neq i$, produced by a spiking of the neuron $i$.

\noindent The autonomous differential equation, verifying some very wide assumptions, governs the system during the inter-spike intervals  of time. It leads to a Poincar\'{e} map (Sotomayor [1979]), which is contractive, as we prove in Theorem \ref{teoremaContraccion} (see also Budelli \em et al.\em [1996]).

\noindent Finally, when synaptic inhibitory coupling is added, it produces asymptotic cyclic attractors in the phase space, whose existence we prove in Theorem \ref{teoremaPrincipal}). Nevertheless, the discontinuities of the synapsis  also generate chaotic orbits, whose existence we prove in Theorem \ref{teoremaCnoVacio}.  While the limit cycles exhibit mostly  asynchronous spikes, but in periodic patterns,  they could seem irregular to an experimental observer, if their periods are very large, and be virtually mixed with the real chaotic orbits.

 \noindent The  model we study includes the networks of coupled  leaky integrators and the relaxation oscillators. The mathematical approach, under particular hypothesis, on  the integrate and fire and the relaxation oscillator models of many neurons was early analyzed by Mirollo \& Strogatz [1990] for homogeneous networks of excitatory neurons. Later  Abbott \& van Vreeswijk [1993] generalized the results for non homogeneous systems, also with delay times and excitatory interactions. For strongly coupled networks of inhibitory or excitatory neurons Bessloff \& Coombes [2000] find synchronous and asynchronous behavior.

  \noindent In the abstract dynamical system
we define the complementary sets $C $ and $S$, of chaotic and stable future orbits respectively (Section 3).

\noindent The qualitative rigorous analysis, discrete-sizing the dynamics by means of a Poincar\'{e} map (see Sotomayor [1979]), was  applied to homogeneous networks by Mirollo \& Strogatz [1990]. Also Budelli \em et al. \em [1991] and   Catsigeras \& Budelli [1992] studied the asymptotic future behavior of inhomogeneous two neurons networks, with a Poincar\'{e}  map in a codimension one section.

\noindent On the other hand discrete-sizings on time,  are  different mathematical approaches and can fit better with computer experiments, in which both time and space are discrete-sized. Cessac [2008] and C\'{e}ssac \& Vi\'{e}ville [2008] obtain similar results to those in this paper, using different mathematical models. They also study discrete discontinuous neural networks. Their model, in spite of being  similar in the evolution given by equation (\ref{ecuacionDiferencial}),  is intrinsically different to consider the spike instants as predetermined by the observer.

\noindent In this paper we  prove that the  Poincar\'{e} map  is discontinuous due to the inhibitory synapsis.
Although it is piecewise contractive, we prove that its discontinuities  play the role of chaos generators: due to them \em the chaotic set $C$ is never empty. \em

\noindent The abstract dynamical systems of this model are classified in three types according to their attractors a.e.. Calling $m$ to the  Lebesgue probability measure, we define the \em  a.e.  chaotic systems \em  (for which $C$ has full measure: $m(C)= 1, \; m(S)= 0$); the \em a.e. stable systems \em (for which $m(C)= 0, \; m(S)= 1$); and \em the combined systems \em (for which $m(C)>0, \; m(S) >0$).

\noindent Generically, under some additional hypothesis,  \em the systems are a.e. stable, and the chaotic set appears as a bifurcation \em among different stable systems. This kind of bifurcations were studied by Catsigeras, Rovella \& Budelli [2008],
proving that they exhibit a chaotic Cantor set attractor.

 \noindent In Theorem \ref{teoremaPrincipal} of this paper we  prove that \em the a.e. stable systems exhibit limit cycles attracting Lebesgue almost all the orbits, \em and due to them they are a.e. periodic. Nevertheless the period of the attracting limit cycles, and the number of them, are mostly determined by the relation among the inhibitory synaptic interactions. They can a priori be very large, and have few relation  with the intrinsic individual periods of the integrate and fire neurons of the network. This last result fits with the  numerical experiments in networks of $n$  different integrate and fire neurons. Postnova \em et al.\em [2007] analyzed the dynamics of computer simulated large networks, and obtained that the system can be driven through different synchronization states, but they are  significantly different from the original periodic behavior  of the individual cells.

 \vspace{.15cm}

\noindent {\bf 1. DEFINITIONS AND STATEMENTS.}

\noindent The phase space is a bounded cube $Q = [-\theta, \theta]^n$ of the vectorial space $\mathbb{R}^n$, whose points are the $n$-uples $V = (V_1, \ldots, V_n)$  with $-\theta \leq V_i \leq \theta, \; \forall \, 1 \leq i \leq n$. We denote $$I= \{1,2,\ldots, n\}$$ to the set of neurons of the network. The variable $V_i$ is the potential states of the $i-$th. neuron of the network, and it  evolves on time $t$.

\noindent In this model the synapses and the spikes of the neurons, can be produced at any instant $t_i$: when the potential $V_i$ of at least one of the neurons (say the $i $-th neuron), reaches $\theta$ (the fixed threshold level).  $t_i $ is not previously specified: it is the solution of the implicit function $\Phi_i^t  = \theta$. We denote as $\Phi^t_i$ to the $i$-th. component of the vectorial flux solution $\Phi^t$  of an autonomous  differential equation \begin{equation} \label{ecuacionDiferencial} dV/dt = F(V) , \;  \forall \, V \in Q \subset \mathbb{R}^n, \;  F:Q \mapsto \mathbb{R}^n, \; F \in {{\cal C}^1} \end{equation}

\noindent This vectorial differential equation  determines the potential $V_i$ of each neuron $i$ as a function on  time $t$, only during  the inter-spike regime.
The solution $\Phi^t(V^0)$  depends on the initial state $V^0= (V^0_1, \ldots, V^0_n)$. We  restrict to the case in which the system is composed by $n$ independent differential equations (i.e. $F_i$ depends only of $V_i$). Nevertheless the open conclusions we obtain (for instance the result in Theorem \ref{teoremaContraccion}) holds  even if the neurons slightly interact during the inter-spike intervals (i.e. $\partial F_i / \partial V_j \approx  0 $ if $i \neq j$). We assume that
\begin{equation}\label{ecuationFsubiPositiva} F_i >0 , \; \;\; \partial F_i/ \partial V_i <0 \; \; \; \forall \; i \in I,\; \; \; \partial F_i / \partial V_j = 0  \; \forall \, j \neq i \end{equation}
and   $F_i$ and $\partial F_i / \partial V_i$  \em are bounded away from zero. \em The dynamics given by (\ref{ecuacionDiferencial}) and (\ref{ecuationFsubiPositiva})  is that of  general   pacemaker   neurons, with strictly increasing potentials in dissipative regime, providing negative concavity of the free evolution on time, of the potential of each neuron in the inter-spike time intervals. This property leads to a contractive Poincar\'{e} map  (Theorem \ref{teoremaContraccion}) and produce, when synaptic inhibitory coupling is added,  asymptotic cyclic attractors in the phase space (see Theorem \ref{teoremaPrincipal}), but also chaotic orbits  (see Theorem \ref{teoremaCnoVacio}).  The limit cycles exhibit mostly  asynchronous spikes, but in periodic patterns.

 \noindent The  model verifying (\ref{ecuacionDiferencial}) and (\ref{ecuationFsubiPositiva}) includes the networks of coupled  leaky integrators and the relaxation oscillators, for instance if $F_i = -\gamma_i (V_i - \beta_i)$, with constant $0<\gamma_i, \beta_i$.

 \noindent From an initial state  $V^0$,  the next spiking instant is   \begin{equation}\label{ecuacionTiempoSpike}\overline {t}(V^0) = \min _{1 \leq i \leq n} t_i(V^0)\end{equation} where $t_i(V^0)$ is the solution of the implicit function \begin{equation}\label{ecuacionTiempoTsubi}\Phi_i^t(V^0)= \theta  \; \; \Leftrightarrow \; \; t = t_i(V^0) \end{equation} At instant $\overline {t}$ maybe more than one neuron reaches the threshold level $\theta$ simultaneously (although this occurs with zero Lebesgue probability in the initial state, if the system is inhibitory).  We call $J(V^0)$ to that set of neurons that spike first, from the initial state $V^0$, i.e.
 \begin{equation} \label{ecuacionJ} J(V^0) = \{i \in I: t_i (V^0) = \overline t(V^0)\} \end{equation} Almost all initial states are such that $J(V^0) = \{i\}$ for a single neuron $i \in I$. The set $J$ of neurons gives  spikes at time $\overline t$, i.e. their potentials  reset to  zero, and from the new initial state $W^0$ it restart the evolution according to the differential equation (\ref{ecuacionDiferencial}). Precisely \begin{equation} \label{ecuationSynapsis1} V_i(\overline t^-) = \theta \; \Rightarrow \; \; V_i(\overline t^+) = W_i(0)= 0 \; \; \forall i \in J\end{equation}

 \noindent When spiking, each of the neurons $i\in J$, which are supposed to be inhibitory, produces a sudden synaptic current through its connections to the other neurons $j \not \in J$. This synapsis produces a sudden change, of amplitude $-H_{ij}<0$, in the potential $V_j$. Precisely: \begin{equation}
 \label{ecuationSynapsis}
 V_j(\overline t^-) < \theta, \; \; \; \; \;   V_j(\overline t^+) = W_j(0)= V_j(\overline t^-) - \sum _{i \in J} H_{ij} \; \; \forall j \not \in J\end{equation}

\vspace{.15cm}

\noindent {\bf 2. THE POINCAR\'{E} SECTION AND THE FIRST RETURN MAP. }

\noindent In this section the dynamics of the first return map $\rho: B \mapsto B$ is  adequately defined  in a Poincar\'{e} section $B$ (see the definition of Poincar\'{e} section and map in Sotomayor [1979]). The section $B$ will be diffeomorphic to the union of a finite number of $n-1$ dimensional balls  transverse to the flux.

\noindent In the compact phase space $Q =  \{V \in \mathbb{R}^n: \; -\theta \leq  V_i  \leq \theta \; \forall \, i \in I\}$ we take the following $n-1$-dimensional section:
\begin{equation}
\label{equationDefiBtechoSubk}
B = \bigcup_{k=1}^n \widehat B_k \; \;  \mbox { where } \; \widehat B_k = \{V \in Q: V_k = 0\}\end{equation}
\noindent Its topology is that induced by $\widehat B_k \subset \mathbb{R}^{n-1}$.

\noindent We assert that \em $B$ is transversal to the flux $\Phi^t$.  \em

\noindent {\em Proof: } It is enough to prove that  $\widehat B_k$ is transversal to the flux for all $k=1, \ldots, n$. This last assertion is deduced from the property $d\Phi_i^t/dt = F_i(\Phi^t) >0$, due to    (\ref{ecuationFsubiPositiva}), as follows: the $n$ components of the vector $d \Phi^t/dt $, which is tangent to the flux, are positive, in particular its $k-$th. component. On the other hand, the manifold $\widehat B_k$ is $n-1$ dimensional, defined by the equation $V_k = 0$. Therefore its tangent subspace $S_k \sim \mathbb{R}^{n-1}$ is formed by all the vectors in $\mathbb{R}^n$ such that have null their $k$-th. component. Then we deduce:
$$S_k \oplus \left [ \frac{d \Phi^t}{dt} \right] = \mathbb{R}^n$$
where $\left[ \vec U \right]$ denotes the subspace generated by the vector $\vec U \in \mathbb{R}^n. \; \; \Box$

\noindent Due to the definition of the mathematical spike in equation (\ref{ecuationSynapsis1}), from any initial state $V^0 \in Q$ the system arrives to the Poincar\'{e} section $B$ in a finite time $\overline t(V^0)$ given by equations (\ref{ecuacionTiempoSpike}) and (\ref{ecuacionTiempoTsubi}). So in particular for $V^0 \in B$, we have defined the first return Poincar\'{e} map $\rho$: \begin{equation} \label{definicionPoincareMap}
\rho: B \mapsto B, \; \; \; \; \rho (V^0) = W^0 = V ^{\overline t^+(V^0)}
\end{equation}
\noindent The sequence of inter-spike intervals (ISI) is given by the evolution of the system through the iterates $\rho^p, p= 1,2, \ldots $ of the return map $\rho$. The ISI sequence is $$\overline t (V^0), \overline t (\rho(V^0)), \overline t(\rho^2(V^0)), \ldots, \overline t (\rho^p(V^0)), \ldots .$$ If the orbits of the discrete dynamical system by  iteration of $\rho$ are attracted to limit cycles or not,  the ISI sequence will be finally periodic or not. If the dynamics of $\rho$ exhibits chaotic attractors, then the ISI sequence also. Even in the periodic case, its period depends of the map $\rho$, (and also of the initial state if there were many limit cycles). But $\rho$  depends strongly of the matrix of synaptic interactions $H_{ij}$. Therefore, the network composed by the inhibitory coupling of $n$ oscillators, may have a dynamics which widely differs from the behavior of each isolated neuron.

\noindent We iterate the Poincar\'{e} map $\rho(V^0) = W^0$, and after each iterate we  reset the time to consider the solution $\Phi^t(W^0)$ of the equation (\ref{ecuacionDiferencial}), from the new initial state $W^0 \in B$. To simplify the notation we will omit the supraindex ${}^0$: from now  $U, V, W$ denote  points in the Poincar\'{e} section.

\noindent From equalities (\ref{ecuationSynapsis1}), (\ref{ecuationSynapsis}) and (\ref{definicionPoincareMap}) we have the following formula for the Poincar\'{e} map:
\begin{eqnarray}
\label{formulaDeRho1}
\rho_j(V)= \; \; &\Phi_j^{\overline t(V)} (V) - \sum_{i \in J(V)} H_{ij}   &\forall \, j \not \in J(V) \\
\rho_i(V) = \;  0 &\; \; &\forall \, i \in J(V) \nonumber
\end{eqnarray}
 \noindent Denote $\wp(I)$ to the family of all non empty parts of the set $I$ of neurons. Roughly speaking $\wp(I)$ the collection of all possible ``words" of \em different \em neurons, of any length $\geq 1, \; \; \leq n$, without considering the order. For each $J \in \wp(I)$ (say for instance, $J = \{1,2\}$), define: $$B_J = \{V \in B: \; \; J(V) = J\}$$ In the example $J=\{1, 2\}$, the set $B_{\{1,2\}} \subset B $ is composed by all the initial states in the Poincar\'{e} section $B$ for which neurons 1 and 2, and only them, will arrive first and simultaneously to the threshold level.

\begin{Remark}
\label{remark0}
\em
\noindent Consider the difference, not only in notation but in significance, among the sets $B_{\{i\}}$ and
$\widehat B_i$.

\noindent On one hand, $B_{\{i\}}$, for a fixed index $i= 1, 2, \ldots, n$, is the set of all the initial states in the Poincar\'{e} section $B$, from which  neuron $i$, and only neuron $i$, will arrive  first to the threshold level, after the   system's evolution in the inter-spike regime. Then $B_{\{i\}} \cap B_{\{j\}} = \emptyset$ if $i \neq j$. Also $\cup_{i=1}^n B_{\{i\}} \subset \neq B = \cup_{J \in {\cap P}(I)} B_J$.

\noindent On the other hand, according to its definition in (\ref{equationDefiBtechoSubk}), $\widehat B_k$,
for a fixed index $k= 1, 2, \ldots, n$, is the set of all the states in the Poincar\'{e} section $B$, in which the potential level $V_k$ of the neuron $k$ is the reset level zero.  We have $\widehat B_k \cap \widehat B_h = \{V \in B: \; V_h = V_k = 0\} \neq \emptyset \; \; \; \forall \, h \neq k$, and $B = \cup_{k=1}^n \widehat B_k$. Also, due to the definition of the Poincar\'{e} map $\rho$ in the equality after (\ref{formulaDeRho1}), we deduce that $\rho(B_{\{i\}}) \subset \widehat B_i, \; \; $    $\rho(B_{\{i,j\}}) \subset \widehat B_i \cap \widehat B_j$.

\end{Remark}

 \begin{Remark} \label{remark1}
 \em
 \noindent Due to formulae (\ref{formulaDeRho1}) the Poincar\'{e} map $\rho$ \em is continuous in each piece $B_J$, in particular in $B_{\{i\}}$ for all $i \in I$. \em  The set $B_J$, for each $J \in \wp(I)$, is called \em  a continuity piece \em of $\rho$.

 \noindent From (\ref{formulaDeRho1}), due to the transversality of the flux $\Phi^t$ to $\widehat B _i$,  $\rho$ transforms homeomorphically  $B_{\{i\}} \cap \widehat B_{k}$ onto its image in $\widehat B_i$. So it is an \em open \em map. Nevertheless it is not globally injective because each point may have different pre-images in $\widehat B_{k}$ and $\widehat B_{h}$ for $h \neq k$.
\end{Remark}

\noindent {\bf The partition ${\cal P}$ and the set  $\partial {\cal P}$ of discontinuities.}

\noindent By construction, for different subsets $J^1 \neq J^2 \in \wp(I)$:  $B_{J_1} \bigcap B_{J_2} = \emptyset $. Any initial state $V$ in $B$ belongs to one and only one set $B_J$ (exactly to that  $B_J$ such that $J(V) = J$). Therefore:

\noindent \em The collection  $\{B_J\}_{ \; J \in \wp(I)}$ form a partition ${\cal P}$ of the Poincar\'{e} section $B$ into the continuity pieces of the return map $\rho$. \em

 \noindent (\ref{ecuacionTiempoSpike}) and (\ref{ecuacionTiempoTsubi}) implies that the set of $V \in B$ such that $\overline t (V) = t_i (V) < t_j(V) \; \; \forall j \neq i$ is open. In other words it is open the set  $B_{\{i\}}$ (i.e. the set of initial states from where the neuron $i$, and only the neuron $i$, arrives first to the threshold level).  On the other hand $B_J$ has empty interior if $\#J \geq 2$.

 \noindent We  consider only those \em not dead \em neurons (i.e. they do not remain under the threshold level forever) and rename as $I$ to that set of neurons. If there is only one, then there is not interactive network to study. So let us suppose that $\# I \geq 2$
 and then $\{B_{{i}} \neq \emptyset : \; \; \;  {i \in I} \}$ is a collection of at least two open and non empty sets. As the space $B$ is connected, we conclude that  \em the topological frontiers of the pieces of the partition ${\cal P}$ is not empty. \em Define:
$$\partial {\cal P} = B \setminus \bigcup _{i=1}^n B_{\{i\}}= \bigcup _{\#J \geq 2} B_J = \bigcup _{i=1}^n \partial B_{\{i\}}.$$

\vspace{.4cm}

\noindent {\bf Topological properties of the Poincar\'{e} map. }

\noindent We use the following norm  to compute the distances in each $\widehat B _k \subset B$, for $k \in I$:
\begin{equation} \label {definicionNorma}\|V - W\| = \max _{i \in I} |V_i - W_i| \; \; \; \forall \, V, W \in \widehat B_k\end{equation}

\noindent Assume the following generic hypothesis in the parameters of the system:
$$H_{ij} \neq \theta \; \; \; \forall \, i \neq j \in I$$
Define the \em expansivity constant: \em
\begin{equation}
\label{definitionAlpha}\alpha  = \frac{\min_{i \neq j} |\theta - H_{ij}|}{4} >0
\end{equation}
\begin{Lemma} \label{lemacontinuidad1}

\noindent If $V  \in \partial {\cal P}$ then  the Poincar\'{e} map $\rho$ is discontinuous in $V$ and the discontinuity jumps in $V$ are larger than $3 \alpha$. Precisely, there exists a sequence of points $U^m \rightarrow V$ such that $\lim _{m \rightarrow + \infty} \|\rho(U^m) - \rho (V)\| > 3 \alpha$ \em
\end{Lemma}

  \noindent{\em Proof: } If $V \in \partial {\cal P}$ then
 $V \in B_J$ for some $J \in \wp(I)$, $\#J \geq 2$. There exist $i \neq j$ in $J$ and therefore $B_J \subset \partial B_{\{i\}}$. By the  definition of topological frontier,  there exists a sequence of points $U^{m} \rightarrow V$ such that $U^{m} \in B_{\{i\}}$.
$$\rho_j(U^{m})  = \Phi_j^{\overline t(U^{m})}(U^{m}) - H_{ij}  \rightarrow _{m \rightarrow + \infty}
\Phi_j^{\overline t(V)}(V)- H_{ij} $$
The last limit is computed recalling that the solution $\Phi^t$ of the differential equation (\ref{ecuacionDiferencial}) depends continuously on the initial state $V$, and the ISI time $\overline t(V)$, defined as the minimum in equality (\ref{ecuacionTiempoSpike}) of the solutions of the implicit equations (\ref{ecuacionTiempoTsubi}), also depends continuously on $V$, due to the Implicit Function Theorem. (Rey Pastor \& col. [1968]).

\noindent But $j \in J, \; V \in B_J \; \; \Rightarrow j \in J(V) $, or in other terms, $V$ is an initial state for which the neuron $j$ arrives to the threshold level. Then $\Phi_j^{\overline t(V)}(V) = \theta$, and by equalities (\ref{formulaDeRho1}), $j \in J(V) \; \Rightarrow \; \rho_j(V) = 0$. Then:
\begin{eqnarray}
\lim \rho_j(U^m) &= &\theta  - H_{ij}, \; \; \; \rho_j(V) = 0 \nonumber \\
|\lim \rho_j(U^m) - \rho_j(V)| &= & |\theta - H_{ij}|  \geq 4 \alpha > 3\alpha \nonumber
\end{eqnarray}
By the definition of limit, for all $m$ large enough we obtain $\|\rho_j(U^m) - \rho_j(V)\| > 3\alpha.  \; \; \; \; \Box$

\noindent The  Lemma \ref{lemacontinuidad1} explains why the discontinuities of $\rho$ act as  saddle type orbits in continuous systems, to produce expansion and chaos. The points of $\partial{\cal P}$ act on their neighborhoods with an infinite rate of expansion in some directions, given by instantaneous expansivity larger than an uniform number $3\alpha >0$. But they also have directions of contraction, as they are in the frontier of the open pieces $B_{\{i\}}$ of continuity of $\rho$, where the following result holds:

{\begin{Theorem}
\label{teoremaContraccion}
The  Poincar\'{e} map $\rho: B \mapsto B$ is uniformly contractive in each of its continuity pieces $B_{\{i\}} \cap \rho(B)$. Precisely, there exists $0 < \lambda < 1$ and a distance $\dist_c$ such that, if
$V, U \in B_{\{i\}}\cap \rho(B)$ for some $i \in I$, then
$$\dist_c(\rho(V),\rho(U)) \leq \lambda \dist_c(V,U)$$
Besides, there exists $K >0$ such that for all $V, U \in B$:$$\frac{1}{K} \|V-U\| \leq \dist_c(V,U) \leq K \|V-U\|$$
\em We note that the distance  $\dist_c$ is not necessarily the norm $\| \cdot \|$ defined in (\ref{definicionNorma}). Nevertheless, they are equivalent in $\mathbb{R}^{n-1}  $. For simplicity we will omit the  not zero constant factors $1/K$ and $K$   in the forward computations along this work, and use $\| \cdot \|$ as if it were  $\dist _c$. This is not a restriction because, as a consequence of  Theorem \ref{teoremaContraccion}, for any norm there exists an iterate $p_0$ such that $\rho^{p_0}$ is contractive in its continuities pieces: $\|\rho^p (V) - \rho^p(U) \| \leq K \dist _c (\rho^p(V), \rho^p(U)) \leq K \lambda ^p \dist_c(V,U) \leq K^2 \lambda ^p \|V-U\| \leq (1/2) \|V-U\|  $ for all $p$ large enough.
\end{Theorem}

\noindent {\em Proof:}
  Define $0<\epsilon_0 = \min _{h \neq j} H_{hj} $. Apply equality (\ref{formulaDeRho1}) and use, at the end of the following formulae, the Lagrange Mean Value Theorem (Rey Pastor \& col. [1968]), to compute the difference of $\Phi_i^t(V)$, for different values of the time $t$, as the increment $\Delta t$ multiplied by the derivative in an intermediate time $T^{m}$:
 \begin{eqnarray}  &V \in \rho(B_{\{k\}})  \; \; \Rightarrow \; \; V_k = 0, \; \; V_j \leq \theta - \epsilon_0 \; \; \forall \; j \in I & \nonumber \\  &V \in \widehat B_ k = \{V \in B: V_k = 0\}  &\nonumber \\ &V \in B_{\{i\}} \; \Rightarrow \;  \phi_i^{\overline t (V)}(V) = \theta, \; \;   \phi_i^{0}(V) = V_i & \nonumber \\ &\epsilon_0 \leq \theta - V_i = \left . \overline t (V) \, d \Phi_i^t(V)/dt \right |_{t = T^{m}} & \nonumber \\ & \epsilon _ 0 \leq \overline t(V) \, F_i ( \Phi^{T^{m}} (V))   & \nonumber \\ & \Rightarrow \; \; \overline t (V) \geq {\epsilon _0}/({\max _{1 \leq i \leq n} \, \max_{V \in Q}  \,  F_i(V)}) = t_0 >0&  \label{formulaTceroMinimo} \end{eqnarray}

\noindent By the Tubular Flux Theorem (Sotomayor [1979]), there exists a ${\cal C}^1$  bounded diffeomorphism which is a spatial change of variables $\xi: V \mapsto \widehat V $ from $Q \subset \mathbb{R}^n$ onto $\widehat Q \subset \mathbb{R}^n$, such that $\xi|_{B} = id$ and the solutions of the differential equation (\ref{ecuacionDiferencial}) in $Q$ verify $d \widehat V/dt = \vec a $ in $\widehat Q$, where $\vec a \in \mathbb{R}^n$ is a constant vector with positive components. It verifies: $\xi (\phi^t(V)) = \xi (V) + \vec a \cdot t, \; \;
d\xi F (V) = \vec a  \; \; \forall \, V \in Q$.

  \noindent Define in ${\mathbb{R}}^n$ the orthogonal projection $\pi$
 onto the $(n-1)$-dimensional subspace
  $a_1 V_1 + a_2 V_2 + \ldots + a_n V_n = 0$. The flux of the differential equation (\ref{ecuacionDiferencial}), after the change $\xi $ of variables in the space,  is ortogonal to that subspace, and is transversal to $\widehat B_k $. Consider any real function $g$:  \begin{eqnarray} \forall\; \widehat{V},  \widehat{V} + d\widehat{V} \in \mathbb{R}^n:\; && \nonumber \\ \pi(d\widehat{V})=  \pi (d\widehat{V} +  \vec a g (\widehat{V}, d \widehat V)) . && \nonumber
  \\ \forall \, V, \, V + dV, \; \; U  \in \widehat B_k, \mbox{ define } && \nonumber \\ \dist_c(V, V + dV) = \| \pi ( d \xi_{V} dV)\| && \nonumber \\ \dist_c( V,U ) = \int _{0}^1 \|\pi \cdot d \xi _{V  + t (U-V)} \cdot (U-V)  \| \, dt && \label{definicionDistanciaContractiva} \end{eqnarray}
  As $\xi$ is a ${\cal C}^1$ diffeomorphism, its derivative and the derivatives of its inverse, are bounded in the compact set $B$, and so the distance $\dist_c$ defined above verifies the last thesis of this Theorem. It is left to prove that $\rho$ is contractive with this distance in $B_{\{i\}} \cap \widehat B _k$.

  \noindent Let us apply $\rho$  to $V, \, V + dV \in B_{\{i\}} \cap \widehat B_k$. We use the equality (\ref{formulaDeRho1})  with $J(V) = \{i\}$, in which for convenience we agree to define $H_{i i} = \theta$. We shall use the derivation formula of the flux of the differential equation  respect to its initial state: \begin{eqnarray} &{\partial \Phi_j ^t (V_j)}/{ \partial V_j} = \exp \left ( {\int _0 ^t } ({\partial F_j }/{ \partial V_j}) \, (\Phi_j^s (V_j)) \, ds) \right )& \nonumber \\ & \mbox{Define: } -\gamma = \max_{1 \leq j} \max _{V \in Q} \; \; {\partial F_j (V_j)}/{ \partial V_j }\; \; < 0  & \nonumber \end{eqnarray}
 \noindent In what follows  $i$ is fixed. It is the value of the index of the continuity piece $B_{\{i\}}$ given in the hypothesis of this lemma. On the other hand, $j = 1, 2, \ldots, i, \ldots, n$ is the index of the general component $\partial \Phi_j^t/\partial t$ of the tangent vector of the flux $\Phi^t$, and of the general component $ \rho_j (V)$ of the  Poincar\'{e} map whose derivative we are computing. We must compute all their components, so we must include the case in which  $j=i$.

  \noindent The formula (\ref{formulaDeRho1}) and the inequality (\ref{formulaTceroMinimo}) lead to:
  $$\rho(V) - \rho(V + dV )  = d \rho(V) \cdot dV =  \left [ ({\partial \rho_j (V)}/{\partial V_j}) dV_j + ({\partial \rho_j(V) }/{\partial V_i}) d V_i   \right ] _{1 \leq j \leq n} 
  {\partial \rho_j (V)}/{\partial V_j}   $$ $$=\left . ({\partial \Phi_j^{t}(V_j) }/{\partial V_j} ) \right |_{t = \overline t (V)} =  \exp \left ( \int _0 ^{\overline t(V)}  ({\partial F_j }/{ \partial V_j}) (\Phi_j^s (V_j)) \, ds) \right )  \leq  e^{- \gamma t_0}      {\partial \rho_j (V)}/{\partial V_i} \,  = $$ $$ = \left . ({d \Phi_j^{t}(V_j) }/{dt}) \right |_{t = \overline t (V)}  \cdot ({d \overline t (V_i)}/{dV_i})=   g(V) \cdot F_j(\Phi^{\overline t (V)}(V)),    $$
  where $g(V) = {d \overline t (V_i)}/{dV_i}$ is the real function obtained deriving respect to $V_i$ the implicit equation $\theta = \Phi_i^{\overline t (V_i)}(V_i) $. Call $\vec e_j$ to the $j-$th. vector of the canonic base in $\mathbb{R}^n$ and join all the results above:

   $$\pi \cdot  d \xi_{\rho(V)}  \cdot (\rho(V) - \rho(V + dV )  ) =  \pi \cdot  d \xi _{\rho(V)} \cdot d \rho(V) \cdot dV =$$ $$= \pi \cdot d \xi _{\rho(V)}  \cdot \left (\sum_{j= 1}^n ({\partial \rho_j(V) }/{ \partial V_j})  \cdot dV_j \vec e_j \right )     + \pi \cdot d \xi _{\rho (V)}  \cdot (g(V) \cdot  F(\phi^{\overline t(V)}(V) ) = $$ $$ = \pi \cdot d \xi _{\rho (V)}  \cdot \left (\sum_{j= 1}^n({\partial \rho_j(V) }/{ \partial V_j})  \cdot dV_j \, \vec e_j \right ) +   g(V) \cdot \pi \cdot d \xi _{\rho(V)}  \cdot F(\phi^{\overline t(V)}(V)) )  = $$ $$ = \pi \cdot d \xi _{\rho (V)} \cdot \left (\sum_{j= 1}^n({\partial \rho_j(V) }/{ \partial V_j})  \cdot dV_j \, \vec e_j \right ) +   \pi ( g(V)  \cdot \vec a )=$$ \begin{equation}=  \pi \cdot d \xi _{\rho(V)} \cdot \left (\sum_{j= 1}^n({\partial \rho_j(V) }/{ \partial V_j} ) \cdot dV_j \, \vec e_j \right )   \label{ecuacionDiferencialDeRho}\end{equation}  Now we define the number $0 <\lambda = e^{-\gamma t_0} <1$ and observe from the computations above that:
  $$0 <\partial \rho(V) / \partial V_j \leq   e^{-\gamma t_0} = \lambda < 1.      $$

\noindent Applying the definition of the differential distance $\dist _c$ in (\ref{definicionDistanciaContractiva}),  and the equality (\ref{ecuacionDiferencialDeRho}), we obtain:
\begin{eqnarray}
&\dist_c(\rho(V), \rho(V + dV) ) = \| \pi (d \xi _{\rho(V)} d \rho_V \cdot dV) \| \leq & \nonumber \\ &  \lambda  \, \|\pi d \xi \cdot dV  \| = \lambda \, \dist_c(V, V + dV) = \lambda \, \| \pi ( d \xi_{V} dV)\|& \nonumber \end{eqnarray}
Integrating by formula (\ref{definicionDistanciaContractiva}) we conclude:
$$\dist _c( \rho (V) , \rho (U) ) \leq \lambda \, \dist_ c (V, U) \; \; \Box$$

\vspace{.1cm}

\noindent {\bf Measure properties of the Poincar\'{e} map.}

\noindent Let $m$ be the $(n-1)$-dimensional  Lebesgue probability measure in the Poincar\'{e} section $B$. Let us prove that

$\ \ \ \ \ \ \ \ \ \ \ \ \ \ \ \ m(\partial{\cal P}) = m  (\bigcup_{\#J\geq 2} B_J ) = 0$.

 \noindent In fact, each $B_J$, if $\#J \geq 2$, is the finite union of  $C^n_2\;$   $n-(\#J)-$dimensional manifolds in $B$, obtained when, for two or more different values of $i \in I$, the respective solutions $t_i (V)$ of  the implicit equations (\ref{ecuacionTiempoTsubi}), are equal.

\noindent Our aim is to study the attractors. The weakest condition required to a set $A \subset B$ to be an attractor is that its basin of attraction has positive Lebesgue measure. So, we may take out the points of the measure zero set $\partial{\cal P}$.

\vspace{.2cm}

\noindent {\bf The set $B'$ of the points with infinite  itinerary.}

\noindent Define

\noindent $  B' = \{V \in B: \rho ^j(V) \not \in \partial{\cal P} \; \; \forall j \geq 0\} \subset B \setminus \partial{\cal P}$

\noindent $B_{\{i\}} \bigcap B_{\{j\}} = \emptyset \mbox{ if } i \neq j \mbox{ and } \bigcup_{i \in I} B_{\{i\}} = B \setminus \partial{\cal P}  \; \Rightarrow $

\noindent $B' = \bigcap _{p= 0}^{+\infty} \rho ^{-p}(B \setminus \partial{\cal P}), \; \; \; \; B \setminus B' = \bigcup _{p = 0} ^{+ \infty} \rho^{-p}(\partial{\cal P}). $

 \noindent The following assertion characterizes the set $B'$:

 \noindent $V \in B' \; \; \Leftrightarrow \; \; \exists  $  a unique sequence $ \{i_p\}_{p \geq 0} $ such that $ \rho^p(V) \in B_{\{i_p\}} \; \forall p \geq 0. $

\noindent This sequence is called \em the itinerary of $V$, \em and it is the infinite sequence of neurons that will spike (in different times), in the order they reach the threshold potential, from the initial state $V$.

 \noindent $\rho$ is continuous  in each $B_{\{i\}}$. Therefore, for all $p \geq 0$ and for all $V \in B'$, the iterate map $\rho^p$ is continuous in $V$. As $B'$ is the numerable intersection of open and dense sets, it is dense.

\begin{Theorem} \label{teoremaBprimaMedidaUno}
The set $B'$ with infinite itinerary  has full Lebesgue measure in $B$,  i.e:
 $m (B \setminus B') = 0$. \end{Theorem}
 \noindent {\em Proof:}   To prove that $m(\rho^{-1}({\partial {\cal P}}))= 0$, once it is known that  $m({\partial{\cal P}}) = 0$,  use the same argument as to prove by induction that for all $p\geq 1:\; m(\rho^{-p}({\partial {\cal P}}))= 0 $, once it is known for $p-1$.

 \noindent In fact, we
 apply Liouville Formula (Sotomayor [1979]) to compute the following Lebesgue integral:
  \begin{eqnarray} 0 = m ({\partial{\cal P}}) &\geq &m ({\rho(\rho^{-1}({\partial{\cal P}}))}) = \nonumber \\
  & = &\int _{V \in \rho^{-1}({\partial{\cal P}})} |\det(D \rho)(V)| \; dm(V) \nonumber \end{eqnarray}
  and conclude that $m((\rho^{-1}({\partial{\cal P}})) \cap\{V \in B; |\det(D \rho)|>0\}) = 0$.
   To prove that $m(\rho^{-1}({\partial{\cal P}}))= 0$ it is enough to prove that
  the Jacobian $det(D \rho)\neq 0$ a.e. in $B$.

  \noindent To apply the Liouville Formula we shall first prove that  $\rho$, given by formulas (\ref{formulaDeRho1}), is differentiable a.e. in the integration set. The first technical problem arises because $\rho $ is not differentiable in all the points: it is neither in the points of discontinuity ${\partial{\cal P}}$ nor in the points of $\bigcup _{h \neq k \in I} (\widehat B_h \cap \widehat B_k)$. In fact, in these last set the Poincar\'{e} section itself fails to be a local differentiable manifold. The set of those exceptional points has zero Lebesgue measure because they are contained in $(n-2)$-dimensional subspaces.

  \noindent The second  technical difficulty  is to check that   $\det D \rho \neq 0 \; a.e. $ The derivative of the Poincar\'{e} map can be computed directly from formulas (\ref{formulaDeRho1}), in each continuity piece $B_{\{i\}}$ intersected with  each $\widehat B_k \setminus \bigcup _{h \neq k}\widehat B_h$. In those sets $B$ is a local differentiable manifold. Computations in (\ref{formulaNombre}) lead to

   \noindent $ \det D \rho(V) = \left [\prod _{j \neq i, \; j \neq k} \;  (\partial \rho_j(V) / \partial V_j) \right ] \cdot (\partial \rho_k(V) / \partial V_i) $ All these factors are computed in equalities (\ref{formulaNombre}) and are not zero due to hypothesis (\ref{ecuationFsubiPositiva}), so  $\det D \rho \neq 0 \; a.e.$ $\Box$

\vspace{.3cm}

\noindent {\bf 3. CHAOTIC AND STABLE SETS.}

\vspace{.1cm}

\noindent We will divide the Poincar\'{e} section $B$ in two complementary sets $S$, formed by stable orbits, and $C = B \setminus S$, the chaotic ones, according to Definition \ref{definicionStableChaoticSets}.
The set $S$ is formed by  stable future orbits under any uniform sufficiently small perturbation that can be added at any step of the iteration of $\rho$. In the set $C = B \setminus S$ there are arbitrarily small perturbations that, if added in some step of the iteration of $\rho$,  drastically change the future orbits and their asymptotic behavior. To have a criteria of chaos, what we call drastic changes in the phase state, we consider the expansivity constant $\alpha>0$ defined in (\ref{definitionAlpha}).

\begin{Definition} \label{definicionStableChaoticSets}
{\bf Stable and chaotic sets.} \em
\noindent $V \in B$ is \em  stable  \em if there exists $\delta >0$ such that
$\forall p \geq 0, \; \;\forall W \in B$,  if $  \|\rho^p(V) - W\| \leq \delta  \; $ then $
\; \| \rho ^k (\rho^p(V)) - \rho ^k (W)\| \leq \alpha  \, \forall k \geq 1.
 $

 \noindent $V \in B$  is \em chaotic \em if it is not stable. The opposite of the definition of stable point holds as follows: $V \in B$ is \em chaotic \em if and only if for all $\delta  >0$ there exists $p \geq 0$, and there exists $W \in B$ such that
$\|\rho^p(V) - W\| \leq \delta  \; $ and $
\; \| \rho ^k (\rho^p(V)) - \rho ^k (W)\| > \alpha  \, \mbox{ for some } k \geq 1.
 $

 \noindent \em  $S$ is the  set of  all the stable orbits   and  $C= B \setminus S$ is the set of all the chaotic orbits. \em

 \noindent It is immediate from the definitions  above that \em  $S$  is forward invariant: \em $\rho(S) \subset S$, and thus, its complement \em  $C$  is backward invariant: $\rho^{-1}(C) \subset C$.
 \end{Definition}

 \noindent Given   $\delta >0$ fixed, we define \em the uniform stable set \em  $S_{\delta } \subset S$   as
 the set of  $V \in B$ such that $\forall p \geq 0, \; \;\forall W \in B, \; \; $ if $  \|\rho^p(V) - W\| \leq \delta  \; $ (with $\delta $ previously fixed), then $
\; \| \rho ^k (\rho^p(V)) - \rho ^k (W)\| \leq \alpha  \, \forall k \geq 1.
 $ From the definitions above observe that $\rho (S_{\delta}) \subset S_{\delta} $ \em  and that any  point in $S$ is in $S_{\delta}$ for some $\delta >0$. \em Therefore, to study the dynamics of the stable points, it is enough to study the dynamics in the sets $S_{\delta}$.

 \begin{Theorem}
 \label{teoremaCnoVacio}

 The set $C$ of chaotic points is never empty. Precisely:
 $C \supset B \setminus B' \supset \partial {\cal P} \neq \emptyset$, where $B'$ is the set of points in the Poincar\'{e} section $B$ with infinite itinerary, and $\partial {\cal P}$ is the set of discontinuities of $\rho$.

 \noindent Even more, if $V \in \partial {\cal P}$ then for all $\delta >0$ there exists $W \in B$ such that
 $\|V - W\| < \delta , \; \; \|\rho(V) - \rho(W)    \| > 3 \alpha$, where $\alpha$ is the expansivity constant of the system.
 \end{Theorem}
 \noindent {\em Proof: }
 From the note at the end of Remark \ref{remark1} and from the definition of the set $B' $ of the points with infinite itinerary, we get $\emptyset \neq \partial {\cal P}  \subset B \setminus B'.$

 \noindent First, let $V \in \partial {\cal P}. $ We shall prove that $V \in C$. From Lemma \ref{lemacontinuidad1} there exists sequence of points $U^{m} \rightarrow V$ such that $\lim \|\rho(U^m ) - \rho (V) \| >  3 \alpha $. Given $\delta  >0$ for all $m\geq 1$ large enough the points $U^m \rightarrow V$ verify $\|U^m - V\| < \delta $. Therefore, taking $W = U^m$ the point $V $ verifies the definition of chaotic point with $p = 0$ and $k=1$.

 \noindent Second, take  $V \in B \setminus B' = \bigcup_{p \geq 0} \rho^{-p}(\partial {\cal P})$. We shall prove that $V \in C$. We know that for some $p \geq 0$ $\rho^p(V) \in \partial{\cal P}. $ In the first step we proved  $\rho^p(V) \in C$, so $V \in \rho^{-p}(C) \subset C$, because $C$ is backward invariant.  $\Box$

\vspace{.1cm}

\noindent {\bf Classification of systems.}

\noindent $\blacktriangleright$ a.e.stable systems:
$m(S) = 1$ and $m(C)= 0$.
In Theorem \ref{teoremaPrincipal} we prove
that the a.e  the limit set is composed only by periodic sinks (limit cycles). Nevertheless the set $C$ of chaotic points in not empty as proved in Theorem \ref{teoremaCnoVacio}.

\noindent $\blacktriangleright$ Chaotic systems: $m(C) = 1$ and $m(S) = 0$.
By definition the limit set $A$ of $C$ is a \em chaotic attractor. \em In Catsigeras \& al. [2008] it is proved that $A$ is a Cantor set attractor.

\noindent $\blacktriangleright$ Partially chaotic and partially periodic systems: $m(C) >0$ and $m(S) >0$. We will show that the points in $S$ are attracted to  limit cycles.

\noindent {\bf Remark: } We did not find ${\cal C}^1$ examples in our model of  systems for which $m(C) >0$ although due to Theorem \ref{teoremaCnoVacio} the set $C$ of chaotic points is never empty.
Nevertheless it is possible (but not immediate) to construct piecewise continuous ${\cal C}^0$ maps $\rho$ in a $n-1\geq 2$ compact ball $B$,  that are uniformly contractive in each of their continuity pieces and such that $C = B, \; \; S = \emptyset $.

\vspace{.1cm}

 \noindent {\bf Some known results.}

\noindent Suppose allowed any small perturbation of the system structure, such that, instead of having the differential vectorial equation (\ref{ecuacionDiferencial}), we have $dV /dt = G (V)$
with $\|G- V\|_{{\cal C}^1} < \epsilon$ for some $\epsilon >0$ sufficiently small. By continuity in the ${\cal C}^1$ topology of the functional space, the first two assumptions (\ref{ecuationFsubiPositiva}) are also verified by $G$ and the third one transforms in $\partial G_i/ \partial V_j \approx 0 \; \; \forall j \neq i $.

\noindent Besides, instead of restricting to a constant matrix $H_0= [H_{ij}]_{i , j \in I}$ to describe the synaptic interactions, allow a matrix $-H(V)$ of negative numbers  $-H( \cdot)< 0$ to be a ${\cal C}^1$ function such that $\|H- H_0\|_{{\cal C}^1} < \varepsilon$.

\noindent The new Poincar\'{e} map $\rho$ will move respect to the old one. The new set of discontinuities $\partial {\cal P}$ of the Poincar\'{e} map  will be (with the Hausdorff distance) near the old one, because its points are those defined by implicit functions, whose equations have ${\cal C}^1$ dependence on $G, H \in {\cal C}^1$. So they are near the old ones if $\|G-F\|_{{\cal C}^1}$ and $\|H - H_0\|_{{\cal C}^1}$ are small enough. We conclude that, if $\varepsilon >0$ is sufficiently small, the new Poincar\'{e} map will be ${\cal C}^1$ near the old one, in each of    its  continuity pieces.

\noindent In this scenario, the following results are known,  leading to the generality  of the property $m(S)=1, m(C)= 0$.

\noindent $\blacktriangleright$ If provided the additional open hypothesis of separation (pairwise disjointness of the closure of the images of the different pieces of continuity for some iterate of $\rho$), then Catsigeras \& al.[2008] proved that generically (open and dense in the ${\cal C}^0$ topology), there exist  a finite number of limit cycles that attracts all the orbits of $B'$. We conclude that $S= B'$, and thus $m(S)= 1$.

\noindent $\blacktriangleright$ If the contractive map $\rho$ is affine in each of its continuity pieces, then Cessac [2008] proved that for generic values of the real parameters, not only in the topological sense (open and dense family of systems), but also in the Lebesgue measure sense (in the space of finite number of the real affinity parameters), there exists at least one and up to  a finite number of limit cycles, which attract all the orbits of $B'$. We conclude that $S= B'$ and $m(S) = 1$.

\noindent $\blacktriangleright$ If the contractive map $\rho$ is not affine, but the matrix of interactions $H$ is constant, Cessac \& Vi\'{e}ville  [2008] proved that either there exists  limit cycles attracting the points of $B'$ or the dynamics has positive entropy.

    \begin{Definition} \em {\bf Omega limit  and limit cycles.}

 \noindent For any $V \in B$, its \em omega-limit  $\omega(V) \in  B$ \em  is the set of limit points of the future orbit of $V$. Precisely:
 $$\omega(V) = \{W \in  B: \exists p_j \rightarrow +\infty \mbox{ such that } \rho^{p_j}(V) \rightarrow W\}$$

 \noindent As the phase space $B$ is compact,  the omega limit set of any point is not empty.  The Poincar\'{e} map $\rho$ is not continuous, so  $\omega(V)$ is not necessarily forward invariant. Nevertheless the omega limit set is the same for all the points in the same orbit, as it is easy to check.

 \noindent A  set $L \in B$ is a \em limit cycle \em (also called \em a periodic sink\em) if it is a single periodic orbit
 and its \em basin attraction $B(L)$  contains an open neighborhood of $L$, \em  being:
 \begin{eqnarray} & L = \{V, \rho(V), \ldots, \rho^{r-1}(V)\}, \; \; \rho^{r}(V) = V \; \;  r \geq 1& \nonumber \\ & B(L) = \{W \in B: \omega (W) = L\}& \nonumber \end{eqnarray}

\end{Definition}

\begin{Theorem}
\label{teoremaPrincipal}
If for some $\delta >0 $ the set $S_{\delta }$ of uniform stable points is not empty then, there exists a finite number $N $ of limit cycles $L_1, \ldots, L_N$ such that the union $\bigcup _{k=1}^N B(L_k)$ of their basins of attractions includes $S_{\delta }$.

\end{Theorem}
\noindent {\bf Remark: } Let us show the main consequence of this theorem: if the set $S$ of stable points has positive measure,  it is not empty and, because all  its points are included in  uniformly stable sets $S_{\delta}$, they all are attracted to  limit cycles. By definition of limit cycles, the basins of attraction of different limit cycles are pairwise disjoint containing open sets. The topology in $B$ has numerable basis, so we conclude that \em there are at most a numerable quantity of limit cycles attracting all the orbits of $S$. \em If besides $m(S) = 1$, then \em $a.e.$ point in $B$ is attracted to a limit cycle. \em

\noindent We will prove Theorem \ref{teoremaPrincipal} at the end of the next Section.

\vspace{.15cm}

\noindent {\bf 4. THE ATOMIZATION OF THE SPACE.}

\noindent Let us consider a system for which the stable set $S $ is not empty. All the results in this Section will hold under this hypothesis.

\noindent As observed at the end of Definition \ref{definicionStableChaoticSets}: $S \neq \emptyset \; \Rightarrow \; S_{\delta } \neq \emptyset $ for some $\delta >0$, and thus also for all sufficiently small $\delta >0$. Fix $0 < \delta < \alpha $ such that $S_{\delta } \neq \emptyset$, where $\alpha $ is the expansivity constant.
We recall that $S_{\delta } \subset B'$ from Theorem \ref{teoremaCnoVacio}, and $\rho (S_{\delta }) \subset S_{\delta}$. Therefore the image $\rho^{p}(S_{\delta })$ is disjoint with $\partial {\cal P}$ for all $p \geq 0$.  So it is partitioned in $n$ disjoint pieces (may be some of them are empty) when intersecting $S_{\delta }$ with $\{B_{\{i\}}\}_{i \in I}$. In each of this pieces  $\rho$ is continuous, and
$\rho(S_{\delta }) = \bigcup _{i \in I} \, \rho (S_{\delta}) \cap B_{\{i\}}$.
\begin{Definition}  {\bf Atoms of generation p.} \em  \label{definitionAtom}
Given $i_1 \in I$  we call \em atom of generation 1 \em  to the following set:

 $ A_{i_1} = \rho(S_{\delta} \cap B_{\{i_1\}}) $

\noindent There are at least one and at most $n = \#I$ non empty atoms of generation 1.
Besides

$\rho (S_{\delta })  = \bigcup _{i_1 \in I} A_{i_1}$.

\noindent Given $(i_1, i_2) \in I \times I$ we call \em atom of generation 2 \em to the following set:
 $A_{i_1, i_2} = \rho (A_{i_1} \cap B_{\{i_2\}})$.

\noindent There are at least one and at most $n^2 = (\#I)^2$ non empty atoms of generation 2.
Besides $\rho^2 (S_{\delta })  = \bigcup _{(i_1, i_2) \in I^2} A_{i_1, i_2}$.

\noindent By induction, if defined the atoms
$\{A_{i_1, i_2, \ldots, i_p}\}_{(i_1, \ldots, i_p) \in I^p}$

\noindent of generation $p$, we define the \em atoms of generation $p+1$: \em
\begin{eqnarray} & A_{i_1, i_2, \ldots, i_p, i_{p+1}} = \rho (A_{i_1, i_2, \ldots, i_p} \cap B_{\{i_{p+1}\}}) &\nonumber \\ & \mbox{ Then: } \; \rho^p (S_{\delta })  = \bigcup _{(i_1, i_2, \ldots, i_p) \in I^p} A_{i_1, i_2, \ldots, i_p}.     & \nonumber \end{eqnarray}

\noindent There are at least one and at most $n^p = (\#I)^p$ non empty atoms of generation $p$.
Fix $p \geq 1$. Denote ${\cal A}_p$  to the finite collection of  atoms of generation $p$:
${\cal A}_p = \{A_{i_1, i_2, \ldots, i_p}\}_{(i_1, \ldots, i_p) \in I^p}$.
\end{Definition}

\noindent {\bf Remark: } It is easy to check from the Definition above that for all $p \geq 1$, if $A \in {\cal A}_{p+1}$ then there exists $A' \in {\cal A}_p$ such that $A \subset A'$. Precisely $A= A_{i_1, i_2\ldots, i_{p+1}} \subset A_{i_2,\ldots, i_{p+1}} = A'$

\noindent We define \em the diameter $d_p$ of the atomization \em of generation $p$: $d_p = \max _{A \in {\cal A}_p} \diam (A)$ where $\diam(A) = \sup_{V, W} \|V- W\|$ denotes the diameter of the atom $A$ and it is 0 if $A = \emptyset$. We observe that  $d_p$ is not the diameter of the union of the atoms, even if they intersect.
\begin{Lemma}
\label{lemadiametroAtomos} Let $0 <\lambda <1$ the uniform contraction rate of Theorem \em \ref{teoremaContraccion}. \em  There exists a constant $K$ such that for any $p \geq 1$ the diameter $d_p$ of the atomization of generation $p$ verifies:

 $0 \leq d_p \leq K \lambda ^{p-1} , \; \; \; \; \; \lim _{p \rightarrow + \infty} d_p = 0$.
\end{Lemma}

\noindent {\em Proof:} Take $K = \mbox{diam} B $. Let us prove the thesis by induction on $p \geq 1$.

\noindent An atom $A \in {\cal {A}}_1$ of generation 1 is the image by $\rho$ of a set in $B$. Then $diam(A) \leq K$. Suppose that all the atoms of generation $p $ have diameter smaller than $K \lambda^{p-1}$. Let us prove the thesis for $p+1$. An atom $A \in {\cal A}_p, \; p \geq 1$ is the image by $\rho$ of some set in $B$, so $A \subset \rho(B)$. Take an atom $A'\in {\cal A}_{p+1}$. By definition:
 $A' = \rho(A \cap B_{\{i\}})$ for some $A \in {\cal A}_p$ and  some $i \in I$.

 \noindent As $A \cap B_{\{i\}}\subset \rho(B) \cap B_{\{i\}}$ due to Theorem \ref{teoremaContraccion} $\rho$  is contractive there:
 \begin{eqnarray} &\mbox{diam}(A') = \sup_{V^1, V^2 \in A'} \|V^1-V^2\| = & \nonumber \\ & =\sup_{U^1,U^2 \in A \cap B_{\{i\}}}
 \|\rho(U^1) - \rho(U^2)\| \leq & \nonumber \\ & \leq  \sup_{U^1,U^2 \in A \cap B_{\{i\}}}\lambda \, \|U^1-U^2\| = & \nonumber \\ &
 =  \lambda \, \mbox{diam}(A \cap B_{\{i\}}) \leq  \lambda \,  \mbox{diam}(A) \leq & \nonumber \\ & \leq  K \, \lambda \, (\lambda^{p-1}) = K\, \lambda ^p. \; \; \Box  & \nonumber\end{eqnarray}
Recall the distance $\dist(V, A)$ of a point $V \in B$ to a non empty set $A \subset B$: $\dist (V,A) = \sup_{W \in A} \|V-W\|$. Denote $\overline A$ and $\partial A$ to the closure and frontier respectively, of the set $A$. A basic classic result from Topology asserts that, in any compact connected metric space, the following properties hold, if $V \not \in A$:
$\dist (V,A) = \dist (V, \overline A ) = \dist (V, \partial A)=  $
$ \max_{W \in \partial A} \|V- W\| $ $= \|V- U_0\| $ for some $U_0 \in \partial A$.
\begin{Lemma} \label{lemaSdeltaDistanciaP}

\noindent If $V \in S_{\delta}$ then \em $\dist (V , \overline{\partial{\cal P}}) \geq \delta /2$.\em
\end{Lemma}
\noindent {\em Proof:} Due to the metric properties in compact connected spaces recalled above, it is enough to show that $\dist (V, {\partial{\cal P}}) \geq \delta  /2$.
By contradiction suppose that $\inf_{U \in \partial{\cal P}} \| V- U\| < \delta  /2$. Then there exists $U \in \partial{\cal P}$ such that:
 $\ \ \ \ \ \ \ \ \ \ \ \ \ \ \ \ \ \ \|V- U \| < {\delta }/{2}$.

\noindent $V \in B_{\delta}$ and by definition of uniform stability, if  $\|V- W\| < \delta $  then:
\begin{equation}
\label{ecuacion}\|\rho^k(V) - \rho^k(W)\| \leq \alpha  \; \; \; \forall k \geq 1
\end{equation}
\noindent In particular $\| \rho(V) - \rho(U) \| \leq \alpha $.

\noindent On the other hand $U \in {\partial{\cal P}}$. Applying Theorem \ref{teoremaCnoVacio} there exists $W \in B$ such that

$\; \; \ \; \ \ \ \ \ \ \ \|U - W\| < {\delta}/{2}, $ $\; \; \|\rho(U) - \rho(W)\| > 3 \alpha$.

\noindent Joining the inequalities above with the triangular property:
$\|V- W\| \leq \|V-U\| + \|U- W\| <  \delta$,

$\; \; 3 \alpha  < \|\rho(U) - \rho(W) \| \leq  $

$\leq \|\rho(U) - \rho(V)\| + \|\rho(V) - \rho(W)\| \leq $

 $ \leq \alpha + \|\rho(V) - \rho(W)\|$,

$\|V- W\| < \delta , \; \; \;$ $ \;  2 \alpha  \leq \|\rho(V) - \rho(W)\|$.
These last two  inequalities contradict (\ref{ecuacion}).$\; \; \; \Box$
\begin{Lemma} \label{lemma44rotulo}
 If $V \in S_{\delta '} \cap B_{\{i\}} $, then  for all $W\in B$ such that  $\|W - V\| < \delta  /2 $ it is verified $  W \in B_{\{i\}}$.

\end{Lemma}

\noindent {\em Proof: } Consider the compact set $K=  B \setminus B_{\{i\}}$. Then:

$U \in K \; \; \Leftrightarrow \; \; \; U \not \in B_{\{i\}} $.

 \noindent We have $V \not \in K$. Call $d = \dist (V, K) = \inf _{U\not \in B_{\{i\}}} \|V -U\|$.
 As recalled when defining the distance of a point to a set,  there exists  $U_0\in \partial K = \partial U_{\{i\}} \subset  \overline{\partial {\cal P}}$ such that

 $d = \| V - U_0\| = \dist (V, \partial K) \geq \dist (V, \overline{\partial {\cal P}})$.

\noindent Due to Lemma \ref{lemaSdeltaDistanciaP}
$\dist (V, \overline{\partial {\cal P}})\geq \delta  /2.$
In resume we have proved that:
$d \geq {\delta}/{2}, \; \; \; \;   \inf _{U \not \in B_{\{i\}}} \|V - U\| \geq {\delta}/2$
Therefore, if  $\|V - W \| < \delta /2 $    then    $  W \in B_{\{i\}}.\; \; \Box$

\vspace{.15cm}
\noindent {\bf Definition. Indivisibility of the atoms. } An  atom $A \in {\cal A}_p$ is \em indivisible \em if its closure $\overline A \subset B_{\{i\}}$ for some $i \in I$.  From the definition of atom, $\rho(A)$ is a unique atom of generation $p+1$, i.e. $A$ does not divide when applying $\rho$. Besides $\rho$ is continuous in $B_{\{i\}}$ so, it is continuous in $\overline A$. Besides, $A \in \rho(B)$, so applying Theorem \ref{teoremaContraccion}:

 \noindent \em If an atom $A$ is indivisible then $\rho$ is continuous and uniformly contractive in $\overline A$.  \em

\vspace{.15cm}

\noindent {\bf Proof of Theorem \ref{teoremaPrincipal}:}

\noindent By Lemma \ref{lemadiametroAtomos} there exists a generation $p\geq 1$ of atoms  such that the diameter $d_p < \delta /2$. Applying Lemma \ref{lemma44rotulo} and recalling that the diameter of a set is the same diameter of its closure, each atom $A \in {\cal A}_p$ is indivisible.  From the remark at the end of Definition \ref{definitionAtom} the atoms of generation $p+1$ and later are contained in the atoms of generation $p$. Then all of them are also indivisible.

 \noindent Fix some non empty atom $A= A^1 \in {\cal A}_p $. As it is indivisible: $\emptyset \neq \rho(A^1) \in {\cal A}_{p+1}$. But any atom of generation $p+1$ is contained in an atom of generation $p$, so there exists $A^2 \in {\cal A}_p$ such that $\rho(A^1)  \subset A^2$. From the indivisibility,  $\rho$ is continuous in $A^1$, thus: $\rho(\overline A^1) \subset \overline A^2$. Applying the same argument to $A_2$ instead of $A_1$ there exists $A^3 \in {\cal A}_p$ such that $\rho (\overline A^2) \subset \overline A^3$.
We deduce $\rho^2(\overline A^1) \subset \rho(\overline A^2) \subset \overline A^3$ where $A^1, A^2, A^3 \in {\cal A}_p$. The family ${\cal A}_p$ is finite, so there is some first pair of integer numbers $k_0\geq 1, r_0\geq 1$ such that $A^{k_0} = A^{r_0+k_0} \in{\cal A}_p$.
 We conclude that  some  subfamily of ${\cal A}_p$ is related in a finite chain:
 \begin{eqnarray}
 \label{ecuacionCadenaAtomos}
 \rho(\overline A^1) \subset \overline A^2, \ldots , \rho(\overline A^{j})\subset \overline A^{j+1},   \ldots,&& \nonumber \\
  \ldots  \rho(\overline A^{k_0}) \subset  \overline A^{k_0 +1},  \ldots,&& \nonumber \\
 \ldots  \rho(\overline A^{k_0 +r_0-1}) \subset \overline  A^{k_0 + r_0} = \overline A^{k_0}
 && \\ \rho^{r_0}(\overline A^{k_0}) \subset \overline A^{k_0} && \nonumber\end{eqnarray}
 All these atoms are indivisible by construction, so for each  $\overline A^{j}, \; \;  \j = 1, 2  \ldots, \leq r_0-1$ in the finite chain  (\ref{ecuacionCadenaAtomos}), the Poincar\'{e} map
 $\rho: \overline{A^j} \mapsto \overline{A^{j+1}}$ is continuous and contractive, with a uniform contraction rate $0 < \lambda < 1$. In resume we have $\rho^{r_0}: \overline A^{k_0} \mapsto \overline A^{k_0} $ and

 \noindent $\|\rho^{r_0}(V) - \rho^{r_0}(W)\| \leq  \lambda^{r_0} \|V- W\|$ for all $ V, W \in \overline A^{k_0}$.

\noindent The Banach Fixed Point Theorem (Lages Lima [1970]), states that in any compact metric space $M$, given a uniformly contractive map $f$ such that $f(M) \subset M$, there exists and is  unique a point $x_1 \in M$  fixed by $f$: $f(x_1) = x_1$. Besides all orbits by future iterates of $f$ have  limit $\{x_1\}$, i.e. $\lim _{k \rightarrow + \infty} f^k(x) = x_1, \; \; \forall x \in M$. We conclude that $\rho^{r_0}$ in the compact metric space $\overline{A^{k_0}} $ has a fixed point $V^1  $ such that:

$V^1 \in \overline{A^{k_0}},$ \ \ $\rho^{r_0}(V^1) = V^1, $

$\lim _{k \rightarrow + \infty} \rho^{  k\,r_0 } (V) = V_1$  $\; \forall V \in \overline{A^{k_0}}$

  \noindent   Therefore $V^1$ is periodic by $\rho$ of period $r_0$, and  its orbit $L_1= \{V^1, \rho(V^1), \ldots, \rho^{r_0 -1}V^1\}$ is the omega limit of all the points in  $\bigcup_{j=1}^{j= k_0 + r_0 } \overline{A}^j$, in particular of those in $\overline A^1$.

\noindent  From the definition of atom and from $\rho(S_{\delta }) \subset S_{\delta}$ each atom is contained in $S_{\delta}$. By Lemma \ref{lemma44rotulo} and using that $\rho$ is contractive in each of its continuity pieces intersected with $\rho(B)$, we deduce that all the points $W$ such that  $\dist (\rho(W), \bigcup_{j=k_0 +1} ^{k_0 + r_0}  \overline A^{j}) < \delta/2$  have $\omega (W) = L_1$. Such points $W$ form an open set ${\cal N} \supset \bigcup_{j=k_0} ^{k_0 + r_0 -1}  \overline A^{j} $. Besides $L_1 \subset  \bigcup_{j=k_0} ^{k_0 + r_0 -1}  \overline A^{j} $. Then ${\cal N}$ is a neighborhood of $L_1$. Therefore
the basin of attraction $B(L_1)$ of the periodic orbit $L_1$ contains a neighborhood of itself, verifying the definition of  limit cycle.

\noindent In particular,  $B(L_1) \supset \overline A^1 \in {\cal A}_p$.
The construction above can be done choosing any first atom $A_1 \in {\cal A}_p$.
We deduce that for each $A \in {\cal A}_p$ there exists a limit cycle $L \subset B'$ such that $\overline A \subset B(L)$.

 \noindent Define
 ${\cal L}= \{L \subset B: L$ is a limit cycle,   $A \subset  B(L) $ for some $ A \in {\cal A}_p \}$. It is a finite  and not empty collection because the family ${\cal A}_p$ of the atoms of generation $p$ is  finite and not empty. By construction the union of the basin of attractions of all the limit cycles
  in $\cal L$ contains all the atoms of generation $p$. Therefore:

 \noindent $\rho^p(B_{\delta } )= \bigcup_{A \in {\cal A}_p} A \subset \bigcup_{L \in {\cal L}} B(L). $
 For all $U $ the limit set $\omega (U) = \omega (\rho^p(U))$. Then $B_{\delta } \subset  \bigcup_{L \in {\cal L}} B(L). \; \; \Box$

\vspace{.1cm}
\noindent {\bf 5. CONCLUSIONS} \smallskip

\noindent We described and analyzed a general mathematical model of a network of $n$ inhibitory pacemaker neurons interacting by synapsis without delay. The size of the network must be finite but as large as wanted. We found  discontinuities in the dynamical system due to the synaptic coupling  and proved that they imply the existence of chaotic orbits. The Poincar\'{e} return map $\rho $ to the $(n-1)$-dimensional section $B$ in the phase space, results from considering the state of the system immediately after each spike. We proved topological and measure properties of $\rho$ as mathematical tools to obtain the dynamical results. We classified the systems according to the measure of the chaotic set and proved that, even being this set never empty, if it has not full measure  there exist stable points and thus,  limit cycles. It is  unknown if there exist ${\cal C}^1$ systems of this model exhibiting a set of chaotic orbits with  positive measure, although under some other additional hypothesis, there are known results about the genericity of the systems exhibiting only limit cycles.

\noindent

\smallskip

\noindent {\bf Acknowledgments} \smallskip

\noindent We thank the project PDT 54/001 of Clemente Estable found (Uruguay),
 the University of Valpara\'{\i}so (Chile) and the University of Marburg (Germany) for partial support, Profs. Pierre Guiraud, Hans Braun and Ruben Budelli for their suggestions, and MEDYFINOL organizing \& scientific comitee for its invitation.  

\smallskip

\noindent {\bf REFERENCES} \smallskip

\noindent{Abbott, L.F. \& Vreeswijk. C [1993]}{`` Assynchronous states in neural networks of pulse-coupled oscillators." }{\it Phys. Rev. E }{\bf 48 } 1483-1490

\smallskip

\noindent{Bessloff, P. \& Coombes, S. [2000]} { `` Dynamics of Strongly Coupled Spiking Neurons."} { \it Neural Computation } {\bf 12} { 91-129.}

\smallskip

 \noindent{Budelli, R. , Catsigeras, E. , Rovella, A. \& G\'{o}mez, L. [1996] }{`` Dynamical behavior of pacemaker neurons networks." }{\it Proc. of the Second Congress of Nonlinear Analysts, Elsevier Science.}

 \smallskip

 \noindent{Budelli, R.,  Torres, J.,    Catsigeras, E. , \&  Enrich, H. [1991] }{`` Two neurons networks I: Integrate and fire pacemakers models." }{\it Biol. Cybern.}{\bf 66},  95-110

\smallskip

 \noindent {Catsigeras, E. \& Budelli, R. [1992] }{`` Limit cycles of a bineuronal  network model. " }{\it Physica D.}{\bf 56}, 235-252

\smallskip

\noindent{Catsigeras, E.,  Rovella, A. \& Budelli, R. [2008]}{`` Contractive piecewise continuous maps modeling networks of inhibitory neurons." }{\it ArXiv [q-bio NC] }{\bf 0805.2695v1}

\smallskip

\noindent{Cessac, B. [2008] }
``A discrete time neural network model with spiking." {\it J. Math. Biol.}  {\bf 54},   311-345

\smallskip

\noindent{Cessac, B. \&  Vi\'{e}ville, T. [2008]  } {``On Dynamics of integrate and fire neural networks with conductance based synapses."} {\it ArXiv [phys-bio-ph]} {\bf 0709.4370v3}

\smallskip

\noindent {Lages Lima, E. [1970]}{``Elementos de Topologia Geral."}
{\it Projeto Euclides, I.M.P.A. , Rio de Janeiro.}

\smallskip

\noindent {Mirollo, R.E. \&  Strogatz, S.H. [1990]  } {``Synchronisation of pulse coupled biological oscillators"} {\it SIAM, J. Appl. Math.} {\bf 50} 1645-1662

\smallskip

\noindent {Postnova, S., Wollweber, B., Voigt, K., Braun, H. A. [2007] }{ ``Neural Impulse Pattern in Bidirectionally Coupled Model of Neurons of Different Dynamics."} {\it Biosystems }{\bf 89 } 135-142

\smallskip

\noindent {Rey Pastor, J., Pi Calleja P.,  Trejo C. [1968] }{``An\'{a}lisis Matem\'{a}tico. Vol. II. C\'{a}lculo infinitesimal de varias variables."} {\it Ed. Kapelusz. Buenos Aires.}

\smallskip

\noindent {Sotomayor, J. [1979]}  {``Li\c{c}\`{o}es de equa\c{c}oes diferenciais ordin\'{a}rias."}
{\it Projeto Euclides, I.M.P.A., Rio de Janeiro.}

\end{document}